\documentclass[11pt,fleqn]{article}
\usepackage{amsmath,amssymb,amsthm} 

\begin{document}

\allowdisplaybreaks

\begin{flushleft}
\LARGE \bf
Lanczos $\tau$-method optimal algorithm in APS\\
 for approximating the mathematical functions
\end{flushleft}

\begin{flushleft} \bf
P.N. Denisenko
\end{flushleft}

\noindent
Kirovograd State Technical University, Kirovograd, Ukraine\\
Email: pnden\_osvita@yahoo.com

\begin{abstract}\noindent
A new procedure is constructed by means of APS in APLAN language.
The procedure solves the initial-value problem for linear
differential equations of order $k$ with polynomial coefficients
and regular singularity in the initialization point in the
interval $[a, b]$ and computes the algebraic polynomial $y_n$
of given order $n$.
A new algorithm of Lanczos $\tau$-method is built for this procedure,
the solution existence $y_n$ of the initial-value problem
proved on this algorithm and also is proved the optimality by precision
of order $k$ derivative of the initial-value problem solution.
\end{abstract}


 \section*{1. ~Introduction}

 {\bf Problem.} To build the procedure with the following parameters:

 {\bf Input:} 1. LDUMK --- Linear differential equation of order $k$
 $$
 LDUMK := ( D[ y ] = 0 ) ; \quad
  D[ y ] = A * y^{(k)} + \cdots  + C * y + G ;           \eqno(1)
  $$
 The coefficients $A , \ldots  , C$ , $G$ of this equation are known polynomials of independent variable $x$.
  The solution $y$ of this equation is a function of variable $x$.
 The zero point is a regular special point of equation (1):

 --- the polynomial $A$ becomes zero in a zero point
 $$
 A( 0 ) = 0 ,                             \eqno(2)
 $$

 --- the equation solution $y$  is a function analytical in the zero point.

 2. Initial conditions in the zero point
 $$
 init\_cond( y , d ) := \{ y( d ) = Y_0 ,  \
 y'( d ) = Y_1 , \
 \ldots ,  \ y^{(l)}( d ) = Y_{l} \};  \ \  d:=0;
 $$
 These conditions are given in a form of Taylor series coefficients partial sum in the zero point
 $$
  T = T_{l}[ y ]  = y(0) + y'(0) *  x + \cdots
 + y^{(l)}(0) / l! * x^{l} ;                       \eqno(3)
 $$

      The initial-value problem for equations (1), (2), (3) meets the following conditions:

 --- the problem has a unique solution

 --- the equation
 $$
 D[ T + V^s[ u ] ] = 0,                                        \eqno(4)
  $$
 where
 $$
 V[u] := \int_0^x u( t ) dt;
  $$
 $$
   s := \max \{ l+1 , k \};                               \eqno(5)
  $$
 $$
   l := deg(T);                                          \eqno(6)
 $$
 $$
 k := ord\_equ(D[y]);                                     \eqno(7)
 $$
 --- equivalent to the problem --- has a recurrence relation system for defining the Taylor coefficients for the function
  $$
 u = c_0 + c_1 * x + \cdots + c_n * x^n + \cdots
 $$
 of the form
  $$
 \{ \ coefTayl( D[  T + V^s[ u  ] ] , i ) = 0 , \quad
  i = 0 , \ldots , m , \ldots \ \} := S_1 \& S_2;
 $$

relating to the function $\ u := y^{(s)} $;

 where $ S_1 $  --- a system of $r$ equations like $0 = 0$;

  $ S_2 $ --- a system of equations without equations like $0 = 0$.

  If the initial-value problem for equations (1), (2), (3) meets these conditions, the solution of this problem is
 $$
 y := solve( LDUMK ,  T) = T + V^s[ u ];              \eqno(8)
 $$

 3. The approximation interval $[a,b]$. Initial point $d \in [a,b]$.

 4. The order $n$ of initial-value problem (1), (2), (3) solution (8) approximation sought within the interval $[a,b]$ is the polynomial
  $$
 y_n = opt. \, \tau\!-\!met( LDUMK , T, [a,b] , n ) =
  c_0 + c_1 * x + \cdots + c_n * x^n                 \eqno(9)
 $$

      {\bf Output.} The polynomial (9). This polynomial meets the initial conditions (3)
 $$
 y_n = T + V^s[ u_p ];              \eqno(10)
 $$
 where
 $$
   u_p = y_n^{(s)};  \ \  p := n - s ;              \eqno(11)
  $$
 --- a polynomial that approximates the equation (4) solution $u=y^ {(s)}$.

 The polynomial $y_n$ (9) derivative of order $k$ within the interval $[a,b]$ approximates the initial-value problem (1), (2), (3) solution (8) derivative of order $k$ optimally by precision
 $$
 \| y^{(k)} - y_n^{(k)} \|_{L_2(a,b;\rho)}
  /  E[ n-k , y^{(k)} , L_2(a,b;\rho) ]
 < Const <  \infty ,                                     \eqno(12)
 $$
 where
 $$
 \| v \|_{L_2(a,b;\rho)}^2 := \int_a^b \ v^2(x) \ \rho(z(x)) \, dx ;
  \ \     \rho(z) := (1 - z^2)^(-1/2);                   \eqno(13)
 $$
 $$
   z(x) := 2*(x-a)/(b-a) - 1;                            \eqno(14)
 $$
  $$
  E[ n-k , y^{(k)} , H ] :=
  \inf_{c_0, \ldots , c_{n-k} }
 \| y^{(k)} - (c_0 + \cdots + c_{n-k}*x^{n-k}) \|_H;        \eqno(15)
$$
 --- the value $y^{(k)}$ of best function approximation by the algebraic polynomials of order $n-k$ in the space $H$

  If the polynomial $y_n$ meets the inequation (12) and  the initial conditions (3), (10), it is optimal by order the joint problem (1), (2), (3) solution $y$ (8) approximation and its derivatives of order $i = 1, 2, \ldots , k$
 --- the following identities take place
 $$
 y^{(i)} - y_n^{(i)} = V^{k-i} [v - v_{n-k}] ; \ \
  i = 0,1, 2, \ldots , k
   $$
 where
$v := y^{(k)}$; $\ v_{n-k} := y_n^{(k)}$
and
 $\| V \|_{C_{[a,b]}} = \max \{|a|,|b|\}$.
\smallskip
\smallskip

 {\bf Relevance of the problem.}
 The equations (1), (2) define the most part of special functions of mathematical physics.
 The order $l$ (6) of polynomial $T$ (3) --- the conditions of initial-value problem for equations (1), (2) --- usually is not equal to $k-1$, where $k$ --- the order (7) of equation (1).
 The Bessel function of order 0 is defined uniquely by the linear differential equation with polynomial coefficients of order 2 and the polynomial (3) of order 0.
 The Frenel integrals are defined uniquely by the linear differential equation with polynomial coefficients of order 3 and the polynomial (3) of order 3.

 Linear differential equations with polynomial coefficients comprise the mathematical modelling apparatus for physical and technical processes [1].
 The computer algebra systems already became the natural mathematical modelling media.
 While analyzing the processes described by equation (1) the computer algebra systems usually perform analytical transformations of equation (1) solution.
 The equation (1) solution usually is not a composition of functions, that can be symbolically transformed by the computer algebra system.
 Thus the computer algebra system transforms the equation (1) solution approximation.
 Usually this approximation is a polynomial.
 The computer algebra systems have the efficient programming means of symbolic polynomial transformations.

 The main criterion for modeling is the model precision.
 According to this criterion the polynomial should approximate the solution optimally by precision order.
 The computer algebra systems have the procedure for computing the Taylor series of  the initial-value problem solution for ordinary differential equations.
 The Taylor series is not the optimal function approximation apparatus --- the order of function $y$ Taylor series remainder term is bigger than the value of this function best approximation
  $$
 \| y - T_n[y] \|_ {C_{[a,b]}}
  /  E[ n , y , C_{[a,b]} ]  = O(q^n), \ \ q > 2
 $$

 Lanczos [1] developed the $\tau$-method of solving the initial-value problem for equation (1).
 Luke applied the Lanczos $\tau$-method for computing the Fourier-Chebyshev coefficients for special functions of mathematical physics.
 These coefficients are the foundation for the procedures of computing the special functions of mathematical physics in mathematical computer software.
 The importance of these procedures initiated the development of new approximative methods of solving the initial-value problem for the equation (1) --- Clenslaw method, Miller method, V.~K.~Dzyadyk a-method [2] and others.

 V.~K.~Dzyadyk a-method [2] solves the the initial-value problem for the equations (1), (3) without the singularity (2) optimally by precision order.
 The analytical methods of solving the initial-value problem for equations (1), (2), (3) optimally by precision are not developed yet.

 \section*{2. ~APLAN-procedure for Lanczos $\tau$-method}

       {\bf Data structure.}

 {\bf Input:} 1. The equation (1) is
 \begin{verbatim}
   LDUMK := ( A * dif( y , k ) + . . . + C * y + G = 0 );
\end{verbatim}
 where $y$ --- an atom, coefficients  $A$ , \dots , $C$, $G$ --- polynomials of equation (atom) $x$.
 These polynomials are separate terms (they are not a multiplication of terms).
 Usually they have the natural form for mathematics and are taken in parenthesis.

       2. Polynomial $T$ (3) has natural form for mathematics
 \begin{verbatim}
             T := d + e * x + . . . + f * x ^ q ;
\end{verbatim}
  The argument $x$ of this polynomial is an atom.

        3. The approximation interval $[a,b]$ description is a list of interval ends.
 \begin{verbatim}
                       interval := (a , b) ;
\end{verbatim}

      {\bf Output.}
       The procedure computes the polynomial $y_n$ (9) with numeric coefficients.
 The form of this polynomial is identical to the form of polynomial $T$.

 {\bf Algorithm 1.}

 1. To compute the equation (1) operator  $D[ y ]$.

 2. To compute the differential operator $D[ y ]$ order $k$ (7).

 3. To compute the order $l$ (6) of the polynomial $T$ (3).

 4. To compute the algorithm parameter $s$ (5).

 5. To compute the order $p$ (11) of the polynomial $u_p$.

 6. To compute the polynomial of order $p$ with symbol coefficients
 $$
 u_p = c_0 + c_1 * x + \cdots  + c_p * x^p;               \eqno(16)
 $$

 7. To compute the polynomial $u_p$ (16) transformation (10).

 8. To compute the transformation of polynomial $y_n$ (16), (10) by operator $D[ y ]$ (1)
  $$
 D[ y_n ] = D[ T + V^s[ u_p ] ] ,                  \eqno(17)
  $$

 9. To compute the zero order of the polynomial $D[ y_n ]$ (17) in the zero point
  $$
  r = deg\_nul(D[  T + V^s[ u_p ]  ]) ;                         \eqno(18)
  $$

 10. To compute the polynomial $D[ y_n ]$ (17) regularization
 $$
  D_0[y_n] = D[  T + V^s[ u_p ]  ] / x^r ;                      \eqno(19)
 $$
  the polynomial $D_0[y_n]$ does not have zeroes in zero point.

 11. To compute the polynomial $D_0[ y_n ]$ (19) order
 $$
 m = deg( D[ T + V^s[ u_p ] ] / x^r) ;                          \eqno(20)
 $$

 12. To compute the auxillary V.~K.~Dzyadyk a-method [2] polynomial of power $m$ (20) and order $p$ (11) (the discrepancy) within the interval $[-1,1]$
 $$
 E_m(x) = E( m , p , x ) = \tau_1 * f_{ p + 1 }( x ) + \cdots
 + \tau_{m-p} * f_m ( x ) ;                         \eqno(21)
  $$
 The discrepancy basis consists of Chebyshev polynomials of first type [2]
  $$
 f_i ( x ) = cheb( i , x ) = cos( i * arccos( x ) ) ;     \eqno(22)
  $$

 13.
To compute the linear transformation (14) of the interval $[ a , b ]$  into the interval $[ -1 , 1 ]$.

 14. To compute the auxillary V.~K.~Dzyadyk a-method [2] polynomial (21) transfer $z(x)$ (14) onto the interval $[a,b]$
  $$
   E_m(z) = E(m,p,z) = subs(x=z,E(m,p,x));   \ \
   z = z(x);                            \eqno(23)
   $$

 15. To compute the V.~K.~Dzyadyk a-method [2] approximation of the left side of the regularized equation (4) --- the sum of polynomials (19), (23)
 $$
 D[  T + V^s[ u_p ]  ]  / x^r + E_m(z)            \eqno(24)
 $$

 16. To compute the system of linear algebraic equations (SLAE)
  $$
 S = \{ \ coefTayl( D[  T + V^s[ u_p  ] ] / x^r + E_m(z) , i ) = 0 , \quad
  i = 0 , \ldots , m \ \} ;                                  \eqno(25)
 $$
 The system variables are the coefficients of the polynomials $u_p$, $E_m$
  $$
 c_0 , \ldots , c_p , \quad \tau_1 , \ldots  , \tau_{m-p}    \eqno(26)
 $$

 17. To compute the SLAE (25) solution --- the coefficients (26) values
 $$
 Coef := solve(S)  := \{\ c_0 = d , \ldots , c_p = e , \
 \tau_1 = f , \ldots  , \tau_{m-p} = g \ \} ;              \eqno(27)
  $$

 18. To compute the polynomial $u_p$ (11) with numeric coefficients
  $$
  \ u_p \  := \ ser( \ Coef  , p \ ) := d +\cdots+ e * x^p;  \eqno(28)
 $$
 The values of the polynomial $u_p$ coefficients (26) are defined by the identities (27).

 19. To compute the polynomial (28) transformation (10).
 This polynomial is the initial-value problem (1), (2), (3) solution (8) approximation sought within the interval $[a,b]$.

 {\bf Algebraic specification of algorithm 1.}

 \begin {verbatim}
    let( LDUMK , Dy = 0 );              /* operator Dy */
    k   := ord_equ( Dy );               /* order  Dy */
    s := deg( T ) + 1;                /* order T + 1 */
    ( k > s ) -> ( s := k );     /* parameter of the method - s */
    p := n - s;                /* power of the polynomial u_p */
    u_p := main_pol(p);   /* u_p with coefficients c(i) */
    y_n := T + n_int( u_p , s );       /* T + V^s[u_p] */
    Dn  := canplf(sub_du(Dy , y_n)); /* polynomial D[y_n] */
    r := deg_nul(canplf(ein_pol(Dn)));  /* zeroes of D[y_n] */
    Dn --> ndiv_x ( Dn , r );             /*  Do[y_n]  */
    m  := deg(canplf(ein_pol(Dn))); /* order Do[y_n] */
    Em := Enl(p, m-p);             /*  E_m for [-1,1]  */
    z --> canplf( -1 + (2/(arg(interval,2) + (-1) *
    arg(interval,1)) * (x + (-1) * arg(interval,1)) );
    Em --> canplf( subs( x = z, Em ));     /*  E_m(z)  */
    Dn --> canplf(Dn + Em);         /*  Do[y_n] + E_m  */
    S := pol_equ(Dn , m);     /*  SLAE - problem approximation */
    Xn := c;    Coef := solve(S);      /* SLAE solution */
    u_p := ser(p , Coef);   /* approximation  of dif(y,s) */
    y_n := T + n_int( u_p , s ); /*  approximation of  y  */
\end{verbatim}

   {\bf The structure of computations results obtained by procedure operators.}

       The procedure solves SLAE (25) relative to the coefficients (26) of form
 \begin{verbatim}
       c(0) , . . . , c(p) , c(p + 1) , . . . , c(m)
\end{verbatim}
 These coefficients are indexed atoms.
 The procedure transforms the polynomials $u_p$ (16) è $E_m$ (21) with these symbol coefficients to obtain the system (25).

         The polynomial (16) with symbol coefficients looks like
 \begin{verbatim}
     u_p := c( 0 ) + c( 1 ) * x + . . . + c( p ) * x ^ p ;
\end{verbatim}
 The result of transforming this polynomial by operators
 $$
 T + n\_int( u\_p , s ) , \ \   canplf(sub\_du(Dy , y\_n)) , \ \
  ndiv\_x ( Dn , r )
 $$
 is a sum of addends like $ c( i ) * x ^ j \$ b$ ,
 where \$ --- the APLAN operation of multiplying the term by constant, constant $b$ is rational and looks like $ rat( p , q ) := p / q$ ;

      The discrepancy (21) with symbol coefficients looks like
 \begin{verbatim}
 E_m := c(p + 1) * cheb(p + 1,x) + . . . + c(m) * cheb(m,x);
\end{verbatim}

 The result of transforming this polynomial by operator \\ $ canplf( subs( x = z, Em ))$ is a sum of addends like
                         $ c( i ) * x ^ j \$ b$ .

       Thus SLAE (25) is a list of equations of form
 \begin{verbatim}
 S := (... ,  c(m) $ f + . . . + c(0) $ e + d = 0 , ... );
\end{verbatim}

       The solution of SLAE (25) by the procedure is a list of identities
 \begin{verbatim}
 Coef := ( c(0) = d , ... , c(p) = f , ... , c(m) = g );
\end{verbatim}

       The procedure transforms the list of identities $Coef$ into the polynomial $u_p$ with numeric natural form coefficients
 \begin{verbatim}
            u_p := d + e * x + . . . + f * x ^ p ;
\end{verbatim}

      {\bf Conclusions from the algebraic specification of the procedure and data structure:}

         1. The procedure has known operators [3].
 These operators perform computations in rational numbers arithmetics.
 The length of numerator and denominator of these numbers is not limited.
 Thus the procedure operators do not add into the computations result the errors of performing the floating point operations.

  2. The procedure has the polynomial complexity by parameter $n$
 $$
 O( n * Q( canplf , n^2 ) ) + O( n^3 ) ,            \eqno(29)
 $$
 where $Q( canplf  , m )$ --- the complexity of transforming by operator
 $canplf$ the polynomial $P$,
 $m$  --- number of addends of the polynomial $canplf( P )$.
 The operator $canplf$ is the internal operator of the APS solver $gr\_solve.exe$.
 It reduces the polynomials into the canonical form for the APS.
 If $P$ is polynomial transformed by the procedure, then
 $canplf( P )$ is a sum of addends like    $ c( i ) * x ^ j \$ b$ .

    {\bf ~Computational experiment with the procedure.}

 The initial-value problem
  $$
 x * y'' - y' + 4 * x^3 * y = 0 , \ \  T \ = \ x^2;              \eqno(30)
 $$
 has the unique solution $y = sin(x)^2$.
 The APLAN-description of initial-value problem (30) and the interval $[-1,1]$ for the procedure implemented above is

 \begin {verbatim}
 process[1] := ( LDUMK :=
 ( x * dif(y , 2) + (-1) * dif(y , 1) + (4 * x^3) * y = 0 );
                T := x ^ 2 ;
         interval := (-1 , 1) ;      ...);
\end{verbatim}

 The results of this initial-value problem transformation by the procedure built with the parameter value $n = 4$ are:

 \begin {verbatim}
 let( LDUMK , Dy = 0 ) ;    Dy := x * dif(y , 2) +
                      (-1) * dif(y , 1) + (4 * x ^ 3) * y ;
 k := ord_equ( Dy ) = 2 ;
 s := deg(T) + 1 = 3 ;
 ( k > s ) -> ( s := k );  s := 3 ;
 p := n - s = 1;
 u_p := main_pol( p ) = c 0 + c 1 * x;
 y_n := T + n_int( u_p , s ) =
 c 1 * x ^ 4 $ rat(1,24) + c 0 * x ^ 3 $ rat(1,6) + x ^ 2 ;
 Dn := canplf( sub_du( Dy , y_n ) ) =
       c 1 * x ^ 7 $ rat(1,6) + c 1 * x ^ 3 $ rat(1,3) +
    c 0 * x^6 $ rat(2,3) + c 0 * x^2 $ rat(1,2) + x^5 $ 4 ;
 r := deg_nul( Dn ) := 2 ;
 Dn --> ndiv_x ( Dn , r ) =
       c 1 * x ^ 5 $ rat(1,6) + c 1 * x $ rat(1,3) +
      c 0 * x ^ 4 $ rat(2,3) + x ^ 3 $ 4 + c 0 $ rat(1,2);
 m  := deg( canplf( ein_pol( Dn ) ) ) = 5 ;
 Em := Enl( p , m - p ) =
       c 5 * x ^ 5 $ 16 + c 5 * x ^ 3 $ -20 + c 5 * x $ 5 +
       c 4 * x ^ 4 $ 8 + c 4 * x ^ 2 $ -8 + c 3 * x ^ 3 $ 4
       + c 3 * x $ -3 + c 2  * x ^ 2 $ 2 + c 4 + c 2 $ -1 ;
 z := x ;
 Em --> canplf( subs( x = z , Em ) ) :=
       c 5 * x ^ 5 $ 16 + c 5 * x ^ 3 $ -20 + c 5 * x $ 5 +
       c 4 * x ^ 4 $ 8 + c 4 * x ^ 2 $ -8 + c 3 * x ^ 3 $ 4
       + c 3 * x $ -3 + c 2  * x ^ 2 $ 2 + c 4 + c 2 $ -1 ;
 Dn --> canplf( Dn + Em ) =
       c 5 * x ^ 5 $ 16 + c 5 * x ^ 3 $ -20 + c 5 * x $ 5 +
       c 4 * x ^ 4 $ 8 + c 4 * x ^ 2 $ -8 + c 3 * x ^ 3 $ 4
       + c 3 * x $ -3 + c 2  * x ^ 2 $ 2 + c 4 + c 2 $ -1 +
       c 1 * x ^ 5 $ rat(1,6) + c 1 * x $ rat(1,3) +
       c 0 * x ^ 4 $ rat(2,3) + x ^ 3 $ 4 + c 0 $ rat(1,2);
 S := pol_equ(Dn , m) = (
       c 4 + c 2 $ -1 + c 0 $ rat(1,2) = 0 ,
       c 5 $ 5 + c 3 $ -3 + c 1 $ rat(1,3) = 0 ,
       c 4 $ -8 + c 2 $ 2 = 0 ,
       c 5 $ -20 + c 3 $ 4 + 4 = 0 ,
       c 4 $ 8 + c 0 $ rat(2,3) = 0 ,
       c 5 $ 16 + c 1 $ rat(1,6) = 0 ) ;
 Coef := solve(S) = ( c 4 = 0 , c 5 = rat(1,14) , c 2 = 0 ,
       c 3 = rat(-9,14) , c 0 = 0 , c 1 = rat(-48,7) ) ;
 u_p := ser(p , Coef) = rat(-48 , 7) * x ;
 y_n := T + n_int( u_p , s ) = x ^ 2 + x ^ 4 $ rat(-1 , 7);
\end{verbatim}

 \section*{3.  ~The solution existence by the procedure}

       {\bf   The Hilbert space $H$.}
 The equation (1) solution is the function analytical throughout the whole complex plane.
 It has the finite number of poles.
 Only the zeroes of the LDEPC (1) coefficients could be the poles of this function.
 Thus it is natural to consider the equation (4) and its approximation (25) in the space of functions analytical in the zero point.
 We will consider the algorithm 1, where the basis meets the following conditions:

 --- the basis elements are the polynomials of power $i$
$$
f_0(x) , \ f_1(x) , \ f_2(x) , \ldots  \ \    deg( f_i(x) ) = i ;   \eqno(31)
 $$

 --- the analytical functions have the Fourier series on the basis (31) within the interval $[a,b]$
  $$
 u  =  c_0 * f_0(z(x)) + \cdots + c_m * f_m(z(x))  + \cdots  ;
 $$
 $$
   v = d_0 * f_0(z(x)) + \cdots  + d_m * f_m(z(x))  + \cdots  ;
 $$
 Thus basis (31) defines the Hilbert space $H$  in the set of analytical functions.
 This space has the scalar product and norm
$$
(  u , v )_H = c_0 * d_0 + \cdots  + c_m * d_m + \cdots  \ \
     || u ||_H^2  = ( u , u )_H                            \eqno(32)
$$

 The basis (22) of V.~K.~Dzyadyk a-method is  is a basis (31) special case.
 It defines the Hilbert space $L_2(-1,1;\rho)$.
 Linear transfer (14) of this basis onto the the interval $[a,b]$ --- $\{ f_m(z(x))\}$ --- defines the space $L_2(a,b;\rho)$.

 The space of algebraic polynomials of order $p$
 $$
 H_p := S_p[H];  \ \
 S_p[c_0 * f_0 + \cdots + c_p * f_p  + \cdots ] :=
  c_0 * f_0 + \cdots + c_p * f_p ;
 $$

--- projecting operator by basis (31),
is a subspace of the space $H$.

       {\bf Theorem 1.}  {\sl
 Let SLAE (25) have the discrepancy basis (31), $p > p_0$ and equation (4) linear regularization operator
 $$
 (D[ T + V^s[ u ] ] / x^r = 0) := ( L[u] + f = 0);             \eqno(33)
  $$
 in the Hilbert space $H$ (32) meet the following conditions:

 --- operator $L$ definitional domain $D( L )$ is dense in the space $H$,

 --- operator $L$ domain of values $R( L )$ is dense in the space $H$,

 --- the operator $L$ transforms the domain $D( L )$ into $R( L)$ biuniquely,

 --- the subspaces $H_p$ --- of algebraic polynomials of order $p$ --- $H_p$ and $L[ H_p ]$ are closed in the $H$ space,

 --- the subspaces sequence $L[ H_p ]$ is dense in $H$ down to the limit;
 $$
 -\!- \  \ \lim_{p \to \infty}  \gamma(p) > 0 , \quad \quad \quad \quad \quad
 \quad \quad  \ \ \gamma(p) :=
 \inf_{  z_p \in L[ H_p ] , \ \|z_p\|_H = 1 } \| S_p[ z_p ] \|_H ;
  $$

     then the SLAE (25) solution $Coef$ (27) exists and it is unique.  }

 {\bf Proof layout. } We proved the equivalence of the SLAE (25) with the basis (31) with equation (33) approximation by projective method with projecting operator by basis (31).
 Thus the theorem 1 is the consequence of main theorem of projectional method convergence.

 \section*{4. ~The solution optimality by the procedure}

    To prove the inequation (12) we transform the integral equation (4) into the linear integral equation of the third type
  $$
   subs(y =  T_{k-1}[y] + V^k[ v ], D[y] = 0 ) :=
 (A*v + \cdots + C*V^k[ v ] + g  = 0) ;                     \eqno(34)
 $$
 The equation (34) solution is the function
  $$
 v = y^{(k)} := V^{s-k}[u] + T^{(k)};
 $$
 The equation (34) is equivalent to the linear integral equation of the second type
 $$
   M[ v ] + g / A = 0 , \ \
   M[v] := v + \cdots + C/A * V^k[ v ] ;      \eqno(35)
  $$
 and the initial conditions
 $$ T_{s-k-1}[v] = T^{(k)}$$

  Similar transformation of the integral equation (4) approximation (24) by V.~K.~Dzyadyk a-method lead to the equation
 $$
   M[ v_{n-k} ] + g / A + x^r / A * E(m,p,z(x)) = 0 ;         \eqno(36)
  $$
 where $ E(m,p,x)$ (21) is the auxillary V.~K.~Dzyadyk a-method [2] polynomial with the basis (31), and the initial conditions
 $$ T_{s-k-1}[v_{n-k}] = T^{(k)}$$

  Equation (36) is the equation (35) approximation by addition of the discrepancy $ x^r / A * E(m,p,z(x))$.
 This discrepancy basis is the transformation of basis (31).
$$
f_0(z(x)) *x^r/A , \ f_1(z(x)) *x^r/A ,
 \ f_2(z(x)) *x^r/A , \ldots        \eqno(37)
 $$

      {\bf Theorem 2.} {\sl
 Let:

 --- linear operator $M[v]$ of the equation (35), $M : H \to H$,
 where $H$ --- the Hilbert space of analytical functions with the scalar product (32) by the basis (31) with the initial conditions
$
  T_{s-k-1}[v] = 0
$
 have reverse operator $ M^{-1}$;

 --- the operator $M^{-1}$  (36) be uniformly limited on the basis (37) functions
 $$
 \| M^{ -1} [ f_i(z(x)) * x^r / A ] \|_H  < W ,
 \ \ \ i = l ,  l+1 , \ldots ;                         \eqno(38)
 $$

 --- SLAE (25) with the basis (31) have the unique solution (27);

 --- for  $  \  n = l+s+1 , l+s+2 , \ldots $  the matrices
 $$
 \{ Q_n \} = \{ ( M^{-1} [ x^r / A * f_{p+i}(z(x)) ] , f_{p+j}(z(x)) )_H
, \ \ \  i , j =  1 , \ldots , m - p \} ,            \eqno(39)
 $$
 where
 $ \ p = n - s , \
 m = \deg( (D[ V^s[ u_p ] + T ]  ) - r ,
 \ r = deg\_nul(  D[ V^s[ u_p ] + T ] ) ;
  $

 have their reverse $ \{ Q_n \}^{-1} $ and the reverse matrices norms
 be uniformly limited by $ n $

 $$
 \| \{ Q_n \}^{ -1} \| := \max_{|C_1|+\cdots+|C_{m-p}|=1 }
 \| \{ Q_n \}^{-1} \{ C_{1} ,  \ldots, C_{m-p}  \}^T \|_{l_2} < Q    \eqno(40)
 $$

  Then assuming $n>l+s$ the inequation (12) is valid and its constant is
 $$
 Const = Q * W
 $$  }

 {\bf Proof layout. }
 The dependence of equation (35) solution $v$ approximation by the equation (36) solution $v_{n-k}$ from the equation (36) discrepancy coefficients
$$
 \| v - v_{n-k} \|_H
  < W * ( |\tau_1 | + \cdots + |\tau_{m-p}| )
 $$
 and the dependence of discrepancy coefficients from Fourier
 coefficients of the equation (35) solution $v$ by the discrepancy basis
  $$
 ( | \tau_1( p )| + \cdots + | \tau_{m-p}( p )| )^2 < Q^2 *
 ( ( v , f_{n-k + 1} )_H^ 2 + \cdots + ( v , f_{n-k+m-p} )_H^ 2 )
  $$

are established.

 Thus the inequation (12) appears directly from these inequations.

     {\bf Example 1.}
 We solved the initial-value problem by the procedure built within the interval $[-1,1]$ with $  n = 4 \ , 6 \ ,  \ldots ,  22$ and computed the polynomials
 $$
 y_n = opt.\ \tau\!-\!method(
 x * y'' - y' + 4 * x^3 * y = 0 , \ y(0) = 0 , \  T  =  x^2
 , \    [-1,1] , \ n ) ;
 $$

  For these polynomials the equation (12) constant $Const$ estimation in the space $C = C_{[-1,1]}$ has the following values
 $$
 \{ \ \| sin''(x^2) - y_n'' \|_C / E[n-2, sin''(x^2), C ] ,
 \ n = 4, 6,\ldots , 22 ; \ \} \ =
 $$
 $$
 \{ \ 1.6 \ , 1.5 \ , 2.1 \ , 1.5 \ , 2 \ , 1.7 \ , 2 \ , 1.7 \ , 2 \ ,
 1.8 \ , 1.9 \ \} ;
    $$

 \section*{5. ~Postamble}

        A new algorithm of solving the initial-value problem for the equations (1), (2), (3) by Lanczos $\tau$-method is built in the article.
 The APLAN language procedure is built according to this algorithm by the programming technology in the APS.
 This procedure solves the initial-value problem for the equations (1), (2), (3) optimally by precision order.
 It illustrates the high efficiency of APS toolkit [3] for creating the new optimal analytical approximative methods of solving the functional equations.


\begin{thebibliography}{3}



 \bibitem{[1]}
 Lanczos C.,
 {\it  Applied Analysis},
 Prendtice Hall, Int. 1956.

\bibitem{[2]}
 Dzjadyk V.K.
 {\it Approximazionnye metody reshenija differenzialnych i
 integralnych uravnenij.}
 Kyiv, Naukova dumka, 1988.

 \bibitem{[3]}
 Denisenko P.N. , Letichevski A.A.
 {\it Algebraicheskoe programmirovanie},
 Kirovograd, KNNPK, 2002.

\end{thebibliography}
 \end{document}